\documentclass{amsart}
\usepackage{amssymb,amsfonts,amsbsy,enumerate,mathrsfs,color}
\usepackage{hyperref}
\usepackage{url}
\usepackage{graphicx}

\newcommand{\G}{\Gamma}
\newcommand{\E}{\overline{\mathcal{E}}}
\newcommand{\MI}{\mathcal{I}}

\theoremstyle{plain}
\newtheorem{theorem}{Theorem}
\newtheorem{corollary}[theorem]{Corollary}

\newtheorem{lemma}[theorem]{Lemma}

\newcommand{\pf}{\noindent{\em Proof: }}

\def\tr{\mathop{\rm tr }\nolimits}

\begin{document}

\title{Regular graphs with maximal energy per vertex}
\author{Edwin R. van Dam}
\address{Department of Econometrics and O.R., Tilburg University,
	P.O. Box 90153, 5000 LE Tilburg, The Netherlands}
\email{Edwin.vanDam@uvt.nl}
\author{Willem H. Haemers}
\address{Department of Econometrics and O.R., Tilburg University,
	P.O. Box 90153, 5000 LE Tilburg, The Netherlands}
\email{Haemers@uvt.nl}
\author{Jack H. Koolen}
\address{School of Mathematical Sciences,
University of Science and Technology of China, Hefei, Anhui 230026 P.R. China} \email{koolen@ustc.edu.cn}
\thanks{This version is published in
Journal of Combinatorial Theory, Series B 107 (2014), 123--131.}

\begin{abstract}
\noindent We study the energy per vertex in regular graphs. For every $k\geq 2$, we give an upper bound for the energy
per vertex of a $k$-regular graph, and show that a graph attains the upper bound if and only if it is the disjoint
union of incidence graphs of projective planes of order $k-1$ or, in case $k=2$, the disjoint union of triangles and
hexagons. For every $k$, we also construct $k$-regular subgraphs of incidence graphs of projective planes for which the
energy per vertex is close to the upper bound. In this way, we show that this upper bound is asymptotically tight.
\end{abstract}




\maketitle


\section{Introduction}

The energy of a graph is the sum of the absolute values of the eigenvalues of its adjacency matrix. This
concept was introduced by Gutman \cite{gutman78} as a way to model the total $\pi$-electron energy of a
molecule. For details and an overview of the results on graph energy, we refer to the recent book by Li, Shi,
and Gutman \cite{energybook} (and the references therein).

Several results on graphs with maximal energy have been obtained. In particular, Koolen and Moulton \cite{km} showed
that a graph on $n$ vertices has energy at most $n(1+\sqrt{n})/2$, and characterized the case of equality. Nikiforov
\cite{Niki} showed that this upper bound is asymptotically tight by constructing graphs on $n$ vertices that have
energy close to the upper bound, for every $n$. Another result, which follows easily from a bound by McClelland
\cite{McC}, is that a graph with $m$ edges has energy at most $2m$, with equality if and only if the graph is the
disjoint union of isolated vertices and $m$ edges (a matching) (see also \cite[Thm.5.2]{energybook}).

In this paper, we consider the (average) energy {\em per vertex} of a graph $\G$, that is,
\begin{equation*}
	\E(\G)=\frac1n \sum_{i=1}^n |\lambda_i|,
\end{equation*}
where $n$ is the number of vertices and $\lambda_1,\lambda_2,\dots,\lambda_n$ are the eigenvalues of $\G$. At
an AIM workshop in 2006, the problem was posed to find upper bounds for the energy per vertex of regular
graphs, and it was conjectured that the incidence graph of a projective plane has maximal energy per vertex;
cf. \cite[Conj.3.11]{workshop}. In this paper, we prove this conjecture
--- among other results.

In Section \ref{sec:mepv} we give an upper bound for the energy per vertex of a $k$-regular graph in terms of $k$, and
show that a graph attains the upper bound if and only if it is the disjoint union of incidence graphs of projective
planes of order $k-1$ or, in case $k=2$, the disjoint union of triangles and hexagons. In order to prove this result,
we reduce the problem to a constrained optimization problem that we solve in Section \ref{sec:opt} using the
Karush-Kuhn-Tucker conditions. Projective planes of order $k-1$ are only known to exist when $k-1$ is a prime power. We
therefore construct, in Section \ref{sec:asy}, $k$-regular subgraphs of incidence graphs of certain elliptic semiplanes
(that are substructures of projective planes) for which the energy per vertex is close to the upper bound, for every
$k$. In this way, we show that our upper bound is asymptotically tight.

\section{Maximal energy per vertex}\label{sec:mepv}

\begin{theorem}\label{thm:epv}
Let $k\geq 2$, and let $\G$ be a $k$-regular graph. Then the energy per vertex of $\G$ is at most
$$\frac{k+(k^2-k)\sqrt{k-1}}{k^2-k+1}$$ with equality if and only if $\G$ is the disjoint union of incidence graphs of
projective planes of order $k-1$ or, in case $k=2$, the disjoint union of triangles and hexagons.
\end{theorem}

\pf First of all, we note that the incidence graph of a projective plane of order $k-1$ (for $k=2$ this is the hexagon)
has spectrum $$\{k^1,\sqrt{k-1}^{k^2-k},-\sqrt{k-1}^{k^2-k}, -k^1\}$$ (see \cite[p. 432]{bcn}), and the triangle has
spectrum $\{2^1,-1^2\}$. Clearly, the disjoint union of several graphs with the same energy per vertex has the same
energy per vertex as the graphs it is built from. The disjoint union of the incidence graphs of projective planes of
order $k-1$ and the disjoint union of triangles and hexagons therefore attain the claimed upper bound for the energy
per vertex. Note also that by considering the disjoint union of several copies of a graph, we may assume without loss
of generality that the number of vertices of the graphs we consider is a multiple of $2(k^2-k+1)$.

Secondly, let $\G$ be a $k$-regular graph with spectrum $\Sigma$. We remark that the bipartite double of $\G$ has
spectrum $\Sigma \cup -\Sigma$ (see \cite[p. 25]{bcn}), and hence it has the same energy per vertex as $\G$. Thus, in
order to prove the claimed upper bound, we may restrict to bipartite graphs, and assume that $\Sigma=-\Sigma$. To show
that this restriction is also possible for the case of equality, we remark that the incidence graph of a projective
plane is not the bipartite double of any graph, except for the hexagon, which is the bipartite double of the triangle.
Indeed, if the bipartite double of $\G$ would be the incidence graph of a projective plane, then $\G$ would have
distinct eigenvalues $k,\sqrt{k-1},$ and $-\sqrt{k-1}$, unless $k=2$, in which case also distinct eigenvalues $2$ and
$-1$ are possible, and $\G$ is a triangle. In the general case, however, $\G$ would be a regular graph with the
property that every pair of vertices has exactly one common neighbor (in other words $\G$ is strongly regular with
$\lambda=\mu=1$), and such graphs do not exist by the Friendship theorem (see \cite[Ch.34]{BOOK}).

In order to find the graphs with maximal energy, we will solve a nonlinear optimization problem that has the
eigenvalues of $\G$ as its variables. Suppose now that $\G$ is a bipartite $k$-regular graph with $n=2m$ vertices and
eigenvalues $\lambda_1 \geq \lambda_2 \geq \cdots \geq \lambda_n$. Because $\lambda_{n-i+1}=-\lambda_{i}$ for all
$i=1,2,\dots,n$, it follows that we only have to consider the first $m$ eigenvalues (which are nonnegative) for the
energy per vertex:

\begin{equation*}
	\E(\G)=\frac1m \sum_{i=1}^m \lambda_i.
\end{equation*}
This is the objective function that we want to maximize. Also in the constraints that we will now formulate, we only
have to consider the first $m$ eigenvalues. The first constraint is that $\tr A^2=nk$, where $A$ is the adjacency
matrix of $\G$. We thus obtain that

\begin{equation*}
	\sum_{i=1}^m \lambda_i^2=mk.
\end{equation*}

The second constraint is obtained by considering $\tr A^4$. This counts the number of walks of length $4$ in the graph,
a number that is at least $nk(2k-1)$: the number of `trivial' walks of length $4$, that is, those not containing a
$4$-cycle. Thus,

\begin{equation*}
	\sum_{i=1}^m \lambda_i^4 \geq mk(2k-1).
\end{equation*}

Together with the constraints that $\lambda_i\leq k$ for $i=1,2,\dots,m$, this gives the optimization problem. We will
solve this problem for fixed $k$ and $m=t(k^2-k+1)$  in the next section. There we will show that the only optimal
solution is $\Sigma^+=\{\lambda_i : i=1,2,\dots,m\}=\{k^t,\sqrt{k-1}^{t(k^2-k)}\}$. Now the only bipartite graphs with
corresponding spectrum $\Sigma=\Sigma^+ \cup -\Sigma^+$ are the disjoint unions of $t$ incidence graphs of projective
planes of order $k-1$ (see~\cite[p.167]{cds82}), which finishes the proof. \qed\\

We remark that an alternative approach, using Cauchy-Schwarz, quickly gives an upper bound of the same order
as our bound. Indeed, if we assume $\G$ to be a connected $k$-regular bipartite graph on $n=2m$ vertices,
then the constraint $\sum_{i=2}^m \lambda_i^2 =(m-k)k$ and Cauchy-Schwarz imply that
$$\E(\G) \leq \frac km+ \frac1m\sqrt{(m-1)(m-k)k}$$
(with equality if and only if $\lambda_2,\lambda_3, \dots, \lambda_m$ are all equal). This upper bound is
increasing in $m$, however. If one could show that $m \leq k^2-k+1$, then that would give an alternative
proof of our result.

\section{Solution of the optimization problem}\label{sec:opt}

In this section, we will solve the optimization problem posed in the proof of Theorem \ref{thm:epv}, where $k$ and
$m=t(k^2-k+1)$ are fixed. In order to do so, we define several functions of ${\bf \lambda}=
(\lambda_1,\lambda_2,\dots,\lambda_m)$. Let $f({\bf \lambda})=\sum_{i=1}^m \lambda_i$, $g=\sum_{i=1}^m \lambda_i^2-mk$,
$h_0({\bf \lambda})=\sum_{i=1}^m \lambda_i^4-mk(2k-1)$, and $h_i({\bf \lambda})=k-\lambda_i$ for $i=1,2,\dots,m$. The
optimization problem under consideration is therefore to maximize $f({\bf \lambda})$ subject to $g({\bf \lambda})=0$,
and $h_i({\bf \lambda})\geq 0$ for $i=0,1,\dots,m$. We first observe that an optimal point $\lambda$ will be nonnegative.

For a feasible point $\lambda$, we let $\MI=\{i=1,2,\dots,m : \lambda_i=k\}$ (the set of indices $i \neq 0$ for which
the constraint on $h_i$ is active). We will now use the Karush-Kuhn-Tucker conditions \cite[p. 246]{per} to examine the
possible optimal points. These conditions state that if a regular feasible point is optimal, then the gradient of $f$
is a linear combination of the gradients of the above defined functions that correspond to the constraints, with
coefficients $0$ for those constraints that are not active. A point is called regular if the gradients corresponding to
the active constraints are linearly independent.

\begin{lemma}\label{lem:three values}
If $\lambda$ is an optimal point such that $h_0({\bf \lambda})>0$, then besides $k$, the entries in $\lambda$ take
at most one value; and if $h_0({\bf \lambda})=0$, then the entries in $\lambda$ take at most two values besides $k$.
\end{lemma}

\pf
We first consider the case that ${\bf \lambda}$ is a nonregular feasible point. If $h_0({\bf
\lambda})>0$, this means that $\nabla g = 2 \lambda$ and $\nabla h_i =-{\bf e}_i$ for $i \in \MI$ are linearly
dependent. Therefore $\lambda_i=0$ for $i \notin \MI$.

If ${\bf \lambda}$ is a nonregular feasible point and $h_0({\bf \lambda})=0$, then $\nabla g = 2 \lambda$,
$\nabla h_0 = 4 \lambda^3$ (defined by $(\lambda^3)_i=\lambda_i^3$), and $\nabla h_i =-{\bf e}_i$ for $i \in \MI$ are
linearly dependent. This implies that there are $c_1$ and $c_2$, not both equal to $0$, such that
$2c_1\lambda_i+4c_2\lambda_i^3=0$ for all $i \notin \MI$. This equation has at most two nonnegative solutions for
$\lambda_i$.

If ${\bf \lambda}$ is a regular and optimal point, then there are $c_1,c_2,$ and $d_i$ for $i=1,2,\dots,m$ such that
$\nabla f +c_1 \nabla g +c_2 \nabla h_0 + \sum_{i=1}^m d_i \nabla h_i=0$, where $c_2=0$ if $h_0({\bf \lambda})>0$ and
$d_i=0$ for $i \notin \MI$. This equation is equivalent to ${\bf j}+2c_1 \lambda +4c_2 \lambda^3 - \sum_{i=1}^m d_i{\bf
e}_j =0$. If $h_0({\bf \lambda})>0$ (and hence $c_2=0$), this reduces to $1+2c_1\lambda_i=0$ for $i\notin \MI$, and
hence $\lambda$ takes at most one value besides $k$. If $h_0({\bf \lambda})=0$, then we find that
$1+2c_1\lambda_i+4c_2\lambda_i^3=0$ for $i\notin \MI$, and this has at most two nonnegative solutions (because of
Descartes' rule of signs, for example). Thus, we proved the lemma. \qed\\

Our next step is to rule out the case that the entries in $\lambda$ take precisely two values besides $k$.

\begin{lemma}\label{lem:two values}
If $\lambda$ is an optimal point, then besides $k$, the entries in $\lambda$ take
at most one value.
\end{lemma}

\pf Suppose on the contrary that this is not the case. Then by Lemma \ref{lem:three values}, we have that $h_0({\bf
\lambda})=0$, and there are --- say --- $m_1$ entries equal to $\theta_1$, $m_2$ entries equal to $\theta_2$, and
$m-m_1-m_2$ entries equal to $k$. For given $k$ and $m$, the optimization problem under consideration can now be
reformulated as
\begin{align*} \label{nlo3}
\text{max.}~~ & m_1\theta_1+m_2\theta_2+(m-m_1-m_2)k \notag \\
\text{s.t.}~~ &  m_1\theta_1^2+m_2\theta_2^2+(m-m_1-m_2)k^2=mk  \notag \\
 &   m_1\theta_1^4+m_2\theta_2^4+(m-m_1-m_2)k^4=mk(2k-1) \notag \\
 & m_1 \geq 0, m_2 \geq 0, m_1+m_2 \leq m, \theta_1 \leq k, \theta_2 \leq k. \notag
\end{align*}
We would like to show that this problem has no optimal solution with $m_1>0$, $m_2>0$, and $\theta_1<\theta_2<k$, so in
the following we only consider such points $x=(m_1,m_2,\theta_1,\theta_2)$. In order to apply again Karush-Kuhn-Tucker,
we let $f'(x)=m_1\theta_1+m_2\theta_2+(m-m_1-m_2)k$, $g'(x)=m_1\theta_1^2+m_2\theta_2^2+(m-m_1-m_2)k^2-mk$,
$h'_0(x)=m_1\theta_1^4+m_2\theta_2^4+(m-m_1-m_2)k^4-mk(2k-1)$, $h'_i(x)=k-\theta_i$ for $i=1,2$, $k'_i(x)=m_i$ for
$i=1,2$, and $k'_3(x)=m-m_1-m_2$.

Now, for a nonregular point $x$ there are $c_1,c_2$, and $c_3$, not all zero, such that $c_1\nabla g'+c_2\nabla
h_0'+c_3\nabla k_3'=0$ (where $c_3=0$ if $m_1+m_2<m$). When we subtract the first two entries in this vector equation
and simplify, we find that $c_1+c_2(\theta_1^2+\theta_2^2)=0$. After dividing the third and fourth entry by $m_1$ and
$m_2$, respectively, and subtracting the resulting equations, we however find (after simplification) that
$c_1+2c_2(\theta_1^2+\theta_1\theta_2+\theta_2^2)=0$. Because $c_2 \neq 0$ (otherwise $c_1=0$, and then also $c_3=0$),
it follows that $(\theta_1+\theta_2)^2=0$, and so $\theta_1=\theta_2=0$ if $x$ is optimal.

For a regular optimal point $x$ (again, satisfying $m_1>0$, $m_2>0$, and $\theta_1<\theta_2<k$), there are $c_1,c_2$,
and $c_3$ such that $\nabla f'=c_1\nabla g'+c_2\nabla h_0'+c_3\nabla k_3'$, where $c_3=0$ if $m_1+m_2<m$. Again,
subtracting the first two entries of this vector equation and simplifying gives that
$1=c_1(\theta_1+\theta_2)+c_2(\theta_1+\theta_2)(\theta_1^2+\theta_2^2)$. Similar as in the nonregular case, we find
from the last two entries of the vector equation that $c_1+2c_2(\theta_1^2+\theta_1\theta_2+\theta_2^2)=0$. The two
obtained equations give that $c_1=2(\theta_1^2+\theta_1\theta_2+\theta_2^2)/(\theta_1+\theta_2)^3$ and
$c_2=-1/(\theta_1+\theta_2)^3$. By substituting these into the third entry of the initial vector equation, and
simplifying, we find that $\theta_1=\theta_2$, and the proof is finished. \qed\\

We note that by allowing the multiplicities to be nonintegral, we have actually proved that optimal points in a
relaxation of the optimization problem take at most one value besides $k$. Thus, we have first generalized the
graph-energy problem in such a way that the eigenvalues could take all possible values at most $k$, and in the previous
step we generalized further by allowing nonintegral multiplicities. Although the hard work has been done now, the proof
is still not finished.

\begin{lemma}\label{lem:the values}
Let $m=t(k^2-k+1)$. If $\lambda$ is an optimal point, then it has $t$ entries equal to $k$, and the remaining entries
equal to $\sqrt{k-1}$.
\end{lemma}

\pf We may now suppose that the entries of an optimal point $\lambda$ take one value --- say --- $\theta$ with
multiplicity $\ell$, besides the value $k$ (with multiplicity $m-\ell$). The optimization problem can thus be
reformulated as
\begin{align*}
\text{max.}~~ & \ell \theta+(m-\ell)k \notag \\
\text{s.t.}~~ & \ell \theta_1^2+(m-\ell)k^2=mk  \notag \\
~ &   \ell \theta_1^4+(m-\ell)k^4 \geq mk(2k-1) \notag \\
~ & 0 \leq \ell \leq m, \theta \leq k. \notag
\end{align*}
This problem is easy enough to be tackled without any sophisticated theory. The first constraint implies that
$\ell=m(k^2-k)/(k^2-\theta^2)$. By substituting this in the second constraint, and simplifying, this constraint reduces
to $\theta \leq \sqrt{k-1}$. Substituting $\ell$ in the objective gives $mk-m(k^2-k)/(k+\theta)$, which is clearly
maximized (subject to the constraints) when $\theta=\sqrt{k-1}$. In this case, the multiplicity of $k$ equals
$m-\ell=t$. \qed

\section{Elliptic semiplanes and an asymptotic result}\label{sec:asy}

A $(k,g)$-cage is a $k$-regular graph with girth $g$ and the smallest possible number of vertices. The incidence graph
of a projective plane of order $q$ is a $(q+1,6)$-cage. In general, it is conceivable that a $(k,6)$-cage has maximal
energy per vertex. For $k-1$ not a prime power only one $(k,6)$-cage is known. The $(7,6)$-cage on $90$ vertices that
was first discovered by Baker \cite{Baker78} is a 3-fold cover of the incidence graph on the points and planes of
$PG(3,2)$; it is the incidence graph of a so-called elliptic semiplane $S(45,7,3)$ (see \cite[p.24, p.210]{bcn}), or
(group) divisible design with parameters $(v,k,\lambda_1,\lambda_2,m,n)=(45,7,0,1,15,3)$. The $(7,6)$-cage has spectrum
\[
\{ 7^1,\ \sqrt{7}^{30},\ 2^{14},\ -2^{14},\ -\sqrt{7}^{30},\ -7^1\}
\]
and hence its energy per vertex is approximately 2.5416.
As a comparison, the upper bound of Theorem~\ref{thm:epv} is approximately 2.5553, and other
graphs close to this bound are the (bipartite double of the) Hoffman-Singleton graph (2.52), the incidence graph of
$AG(2,7)$ minus a pencil (2.4965),the incidence graph of $AG(2,7)$ minus a parallel class (2.4106),
the incidence graph of the biplane on 29 points (2.4003), and the Klein graph (2.3472).

The two examples from affine planes generalize.
An affine plane of order $q$ from which a parallel class is deleted is an elliptic
semiplane $S(q^2,q,q)$.
The incidence graph $\G$ of such an elliptic semiplane has spectrum
\[ \{ q^1, \ \sqrt{q}^{q(q-1)},\ 0^{2(q-1)}, \ -\sqrt{q}^{q(q-1)},\ -q^1 \},
\]
and therefore $\E(\G)=\sqrt{k}-\frac{1}{\sqrt{k}}+\frac{1}{k}$.
An affine plane of order $q$ from which a pencil (one point $x$ together with all lines through $x$) is deleted,
is an elliptic semiplane $S(q^2-1,q,q-1)$ whose incidence graph $\G$ has spectrum
\[ \{ q^1, \ \sqrt{q}^{q^2-q-2},\ 1^{q},\  -1^q,\ -\sqrt{q}^{q^2-q-2},\ -q^1 \},
\]
and hence ${\E}(\G)=\sqrt{k}+\frac{2k}{k^2-1}-\frac{\sqrt{k}}{k-1}$, which is slightly better than the first example.
For $k=11$ the latter formula gives $\E(\G)\approx 3.1683$, whilst the upper bound of Theorem~\ref{thm:epv} is
approximately 3.2329. Both families obtained from affine planes have girth 6, but the second one has two vertices less.

The mentioned families also lead to examples of $k$-regular graphs with large energy per vertex for arbitrary $k$.
Indeed, for both types of elliptic semiplanes one can delete $\ell=q-k$ point classes and $\ell$ line classes such that
the remaining structure is a square $1$-design with block size~$k$. From the elliptic semiplane $S(q^2-1,q,q-1)$ we
thus obtain a $k$-regular incidence graph $\G$ on $2(q^2-1-\ell(q-1))$ vertices. Eigenvalue interlacing
(see~\cite[p.37]{bh12}) gives that except for $\pm k$, $\G$ has eigenvalues $\pm\sqrt{q}$, both with multiplicity at
least $q^2-q-2-2\ell(q-1)$, and at least $2q$ additional eigenvalues which are in absolute value at least 1. Hence
\[
\E(\G)\geq \frac{2q-\ell+\sqrt{q}(q^2-2q\ell+2\ell-q-2)}{(q-1)(q+1-\ell)}.
\]
This leads to the following lower bound.
\begin{theorem}\label{k}
Let $k$ be a positive integer, and let $\ell$ be the smallest nonnegative integer such that $k+\ell$ is a prime power.
Then there exists a $k$-regular graph $\G$ whose energy per vertex $\E(\G)$ satisfies
\[
\E(\G)\geq \frac{2k+\ell+\sqrt{k+\ell}(k^2-\ell^2+\ell-k-2)}{(k+1)(k+\ell-1)}.
\]
\end{theorem}
As a corollary we find that the bound of Theorem~\ref{thm:epv} is asymptotically tight:
\begin{corollary}
If $k$ is large enough, then there exists a $k$-regular graph $\G$ for which
\[
{\E}(\G) \geq \sqrt{k}-k^{1/40}.
\]
\end{corollary}
\pf
It is known that $\ell\leq k^{21/40}$ if $k$ is large enough (see \cite{bhp}),
and for $k\rightarrow\infty$ the $k$-regular graph of Theorem~\ref{k} satisfies
\[\pushQED{\qed}
\E(\G)\geq\sqrt{k}-\frac{\ell}{2\sqrt{k}}- o(1). \qedhere
\popQED\]

\noindent{\em Remark.} Under the Riemann Hypothesis, it can be proved that for large $k$ there is always a prime number
between $k$ and $k+\sqrt{k}$. This would improve the lower bound of the corollary to $\sqrt{k}-\frac12-o(1)$.

We also note that the graphs that attain the Koolen and Moulton \cite{km} bound for the energy are $k$-regular
with $k=(n+\sqrt{n})/2$, and hence the energy per vertex of such graphs is approximately $\sqrt{k/2}$.
Nikiforov's \cite{Niki} related examples are (not necessarily regular) subgraphs of Paley graphs, and the
latter also have energy per vertex approximately $\sqrt{k/2}$.

\section{Conclusion and final remarks}

In this paper, we obtained a bound on the energy per vertex in a $k$-regular graph, thus proving a conjecture posed at
an AIM workshop \cite{workshop} in 2006. The incidence graphs of projective planes of order $k-1$ attain this bound.
For values of $k$ for which no projective plane of order $k-1$ exists, we construct $k$-regular graphs from elliptic
semiplanes for which the energy per vertex is close to the bound, and show in this way that the bound is asymptotically
tight.

We note that on the other extreme, it is relatively easy to show that the energy per vertex of a $k$-regular
graph is at least $1$, with equality if and only if the graph is a disjoint union of copies of the complete
bipartite graph $K_{k,k}$ (indeed: $\frac1m \sum_{i=1}^m \lambda_i \geq \frac1{mk} \sum_{i=1}^m
\lambda_i^2=1$ with equality if and only if the (nonnegative) eigenvalues are $0$ or $k$).

Fiorini and Lazebnik \cite{FL98} showed that the incidence graphs of projective planes are also extremal in
the sense that they have the largest number of $6$-cycles among the bipartite graphs (on $m+m$ vertices)
without $4$-cycles. De Winter, Lazebnik, and Verstra\"{e}te \cite{DwLV08} obtained a similar result for
$8$-cycles.

It would also be interesting to study the energy per vertex for graphs that are not necessarily regular. For
example, if we consider the energy per vertex of trees, then it follows from the fact that a path on $n$
vertices has maximal energy among all trees on $n$ vertices, and an expression for its energy (see
\cite[p.26]{energybook}) that the energy per vertex of trees is less than $4/\pi$, and that this bound is
tight.\\

\noindent {\bf Acknowledgements.} The authors thank the referees for their useful comments. JHK was partially
supported by the 100 talents program of the Chinese government. This work was done while JHK was visiting the
Department of Econometrics and Operations Research of Tilburg University, for which support from NWO is
gratefully acknowledged.


\begin{thebibliography}{99}

\bibitem{BOOK} M. Aigner and G.M. Ziegler, \emph{Proofs from THE BOOK},
    Springer, third edition, 2004.

\bibitem{bhp} R.C. Baker, G. Harman, and J. Pintz, The difference between consecutive primes, II, {\sl Proc. Lond.
    Math. Soc. (3)} 83 (2001), 532--562.

\bibitem{Baker78} R.D. Baker, An elliptic semiplane, {\sl J. Combin. Theory Ser. A} 25 (1978), 193-–195.

\bibitem{bcn} A.E. Brouwer, A.M. Cohen, and A. Neumaier,
    \emph{Distance-Regular Graphs}, Springer-Verlag, Berlin-New York, 1989.

\bibitem{bh12} A.E. Brouwer and W.H. Haemers, \emph{Spectra of Graphs},
    Springer, 2012; available online at \url{http://homepages.cwi.nl/~aeb/math/ipm/}.

\bibitem{workshop} R. Brualdi, L. Hogben, and B. Shader, AIM workshop spectra of families of matrices described by
    graphs, digraphs, and sign patterns - final report: mathematical results, 2007; available online at
    \url{http://aimath.org/pastworkshops/matrixspectrumrep.pdf}.

\bibitem{cds82} D.M. Cvetkovi\'c, M. Doob, and H. Sachs, \emph{Spectra of
    Graphs}, VEB Deutscher Verlag der Wissenschaften, Berlin, 1979.

\bibitem{DwLV08} S. De Winter, F. Lazebnik, and J. Verstra\"{e}te, An extremal characterization of projective planes, {\sl Electron. J. Combin.} 15 (2008), R143.

\bibitem{FL98} G. Fiorini and F. Lazebnik, An extremal characterization of incidence graphs of projective
    planes, {\sl Acta Appl. Math.} 52 (1998), 257--260.

\bibitem{gutman78} I. Gutman, The energy of a graph, {\sl Ber. Math.-Stat. Sekt. Forschungszent. Graz} 103 (1978),
    1-–22.

\bibitem{km} J.H. Koolen and V. Moulton, Maximal energy graphs, {\sl Adv. in Appl. Math.} 26 (2001),
    47-–52.

\bibitem{energybook} X. Li, Y. Shi, I. Gutman, \emph{Graph Energy}, Springer, 2012.

\bibitem{McC} B. McClelland, Properties of the latent roots of a matrix: The estimation of $\pi$-electron energies, {\sl
    J. Chem. Phys.} 54 (1971), 640--643.

\bibitem{Niki} V. Nikiforov, Graphs and matrices with maximal energy, {\sl J. Math. Anal. Appl.} 327 (2007), 735-–738.

\bibitem{per} A.L. Peressini, F.E. Sullivan, and J.J. Uhl Jr., \emph{The Mathematics of Nonlinear
    Programming}, Springer, 1988.

\end{thebibliography}
\end{document}